\documentclass[12pt]{article}
\usepackage{e-jc}
\usepackage{amssymb,amsthm,amsfonts,amsmath}

\usepackage{hyperref}
\usepackage{mathrsfs}

\usepackage[dvips]{graphicx}
\usepackage{wrapfig}

\newcommand{\R}{\mathbb R}
\newcommand{\Y}{\mathbb Y}

\newcommand{\A}{\mathbb A}

\newcommand{\al}{\alpha}

\newcommand{\ka}{\varkappa}
\newcommand{\la}{\lambda}
\newcommand{\La}{\Lambda}
\newcommand{\de}{\delta}

\newcommand{\tth}{\theta}
\newcommand{\si}{\sigma}

\newcommand{\pd}{\partial}

\newcommand{\gr}{\operatorname{gr}}

\newcommand{\up}{\uparrow}
\newcommand{\down}{\downarrow}

\newcommand{\HH}{\mathbf H}
\newcommand{\h}{\mathbf h}
\newcommand{\p}{\mathbf p}

\newtheorem{theorem}{Theorem}[section]
\newtheorem{proposition}[theorem] {Proposition}
\newtheorem{corollary}[theorem]{Corollary}

\theoremstyle{definition}
\newtheorem{definition}[theorem]{Definition}
\newtheorem{remark}[theorem]{Remark}

\numberwithin{equation}{section}

\title{Plancherel averages: \\Remarks on a paper by Stanley}

\author{Grigori Olshanski\thanks{Supported by a grant from the Utrecht University,
by the RFBR grant 08-01-00110, and by the project SFB 701
(Bielefeld University).}\\
\small Institute for Information Transmission Problems\\[-0.8ex]
\small Bolshoy Karetny 19\\[-0.8ex]
\small Moscow 127994, GSP-4, Russia \\[-0.8ex]
\small and\\[-0.8ex]
\small Independent University of Moscow, Russia \\
\small \texttt{olsh2007@gmail.com} }

\date{\dateline{October 1, 2009}{March 10, 2010}}

\begin{document}
\maketitle

\begin{abstract}
Let $M_n$ stand for the Plancherel measure on $\Y_n$, the set of Young diagrams
with $n$ boxes. A recent result of R.~P.~Stanley (arXiv:0807.0383) says that
for certain functions $G$ defined on the set $\Y$ of all Young diagrams, the
average of $G$ with respect to $M_n$ depends on $n$ polynomially. We propose
two other proofs of this result together with a generalization to the Jack
deformation of the Plancherel measure.
\end{abstract}

\section{Introduction}\label{1}

Let $\Y$ denote the set of all integer partitions, which we identify with Young
diagrams. For $\la\in\Y$, denote by $|\la|$ the number of boxes in $\la$ and by
$\dim\la$ the number of standard tableaux of shape $\la$. Let also
$c_1(\la),\dots,c_{|\la|}(\la)$ be the contents of the boxes of $\la$ written
in an arbitrary order (recall that the content of a box is the difference $j-i$
between its column number $j$ and row number $i$).

For each $n=1,2,\dots$, denote by $\Y_n\subset\Y$ the (finite) set of diagrams
with $n$ boxes. The well-known {\it Plancherel measure\/} on $\Y_n$ assigns
weight $(\dim\la)^2/n!$ to a diagram $\la\in\Y_n$. This is a probability
measure. Given a function $F$ on the set $\Y$ of all Young diagrams, let us
define the $n$th {\it Plancherel average\/} of $F$ as
\begin{equation}\label{1A}
\langle F\rangle_n=\sum_{\la\in\Y_n}\frac{(\dim\la)^2}{n!}\,F(\la).
\end{equation}

In the recent paper \cite{St3}, R.~P.~Stanley proves, among other things, the
following result (\cite[Theorem 2.1]{St3}):

\begin{theorem}\label{thm1.A}
Let $\varphi(x_1,x_2,\dots)$ be an arbitrary symmetric function and set
\begin{equation}\label{1B}
G_\varphi(\la)=\varphi(c_1(\la),\dots,c_{|\la|}(\la),0,0,\dots), \qquad
\la\in\Y.
\end{equation}
Then $\langle G_\varphi\rangle_n$ is a polynomial function in $n$.
\end{theorem}

The aim of the present note is to propose two other proofs of this result and a
generalization, which is related to the {\it Jack deformation\/} of the
Plancherel measure.

The first proof relies on a claim concerning the shifted (aka interpolation)
Schur and Jack polynomials, established in \cite{OO1} and \cite{OO2}. Modulo
this claim, the argument is almost trivial.

The second proof is more involved but can be made completely self-contained. In
particular, no information on Jack polynomials is required. The argument is
based on a remarkable idea due to S.~Kerov \cite{Ke2} and some considerations
from my paper \cite{Ol}.

As indicated by R.~P.~Stanley, his paper was motivated by a conjecture in the
paper \cite{Han} by G.-N.~Han (see Conjecture 3.1 in \cite{Han}). Note also
that examples of the Plancherel averages of functions of type \eqref{1B}
appeared in S.~Fujii et al. \cite[Section 3 and Appendix]{Fuj}.

\section{The algebra $\A$ of regular functions on $\Y$}\label{2}

For a Young diagram $\la\in\Y$, denote by $\la_i$ its $i$th row length.
Clearly, $\la_i$ vanishes for $i$ large enough. Thus, $(\la_1,\la_2,\dots)$ is
the partition corresponding to $\la$.

\begin{definition}\label{defn2.A}
Let $u$ be a complex variable. The {\it characteristic function\/} of a
diagram $\la\in\Y$ is
$$
\Phi(u;\la)=\prod_{i=1}^\infty\frac{u+i}{u-\la_i+i}
=\prod_{i=1}^{\ell(\la)}\frac{u+i}{u-\la_i+i}\,,
$$
where $\ell(\la)$ is the number of nonzero rows in $\la$.
\end{definition}

The characteristic function is rational and takes the value 1 at $u=\infty$.
Therefore, it admits the Taylor expansion at $u=\infty$ with respect to the
variable $u^{-1}$. Likewise, such an expansion also exists for $\log
\Phi(u;\la)$.

\begin{definition}\label{defn2.B}
Let $\A$ be the unital $\R$-algebra of functions on $\Y$ generated by the
coefficients of the Taylor expansion at $u=\infty$ of the characteristic
function $\Phi(u;\la)$ (or, equivalently, of $\log\Phi(u;\la)$). We call $\A$
the {\it algebra of regular functions on $\Y$}. (In \cite{KeO} and \cite{IO},
we employed the term {\it polynomial functions\/} on $\Y$.)
\end{definition}

The Taylor expansion of $\log\Phi(u;\la)$ at $u=\infty$ has the form
$$
\log\Phi(u;\la)=\sum_{m=1}^\infty\frac{p^*_m(\la)}m\,u^{-m},
$$
where, by definition,
$$
p^*_m(\la)=\sum_{i=1}^\infty [(\la_i-i)^m-(-i)^m]=\sum_{i=1}^{\ell(\la)}
[(\la_i-i)^m-(-i)^m], \qquad m=1,2,\dots, \quad \la\in\Y.
$$
Thus, the algebra $\A$ is generated by the functions $p^*_1,p^*_2,\dots$. It is
readily verified that these functions are algebraically independent, so that
$\A$ is isomorphic to the algebra of polynomials in the variables
$p^*_1,p^*_2,\dots$. Note that $p^*_1(\la)=|\la|$.

Using the isomorphism between $\A$ and $\R[p^*_1,p^*_2,\dots]$ we define a {\it
filtration\/} in $\A$ by setting $\deg p^*_m(\,\cdot\,)=m$. In more detail, the
$m$th term of the filtration, consisting of elements of degree $\le m$,
$m=1,2,\dots$, is the finite-dimensional subspace $\A^{(m)}\subset\A$ defined
in the following way:
$$
\A^{(0)}=\R1; \quad
\A^{(m)}=\operatorname{span}\{(p^*_1)^{r_1}(p^*_2)^{r_2}\dots\,:\,
1r_1+2r_2+\dots\le m\}.
$$

The regular functions on $\Y$ (that is, elements of $\A$) coincide with the
{\it shifted symmetric functions\/} in the variables $\la_1,\la_2,\dots$ as
defined in \cite[Sect. 1]{OO1}. Thus, we have the canonical isomorphism of
filtered algebras $\A\simeq\La^*$, where $\La^*$ stands for the algebra of
shifted symmetric functions. This also establishes an isomorphism of graded
algebras $\gr\A\simeq\La$, where $\La$ denotes the algebra of symmetric
functions.

For a diagram $\la\in\Y$, denote by $\de(\la)$ the number of its diagonal
boxes, by $\la'$ the transposed diagram,  and set
\begin{equation}\label{2A}
a_i=\la_i-i+\tfrac12, \quad b_i=\la'_i-i+\tfrac12, \qquad i=1,\dots,\de(\la).
\end{equation}
We call the numbers \eqref{2A} the {\it modified Frobenius coordinates\/} of
$\la$ (see \cite[(10)]{VK}).

\begin{proposition}\label{prop2.C}
Equivalently, $\A$ may be defined as the algebra of {\it super-symmetric
functions\/} in the variables $\{a_i\}$ and $\{-b_i\}$.
\end{proposition}

\begin{proof}
See \cite{KeO}. Here I am sketching another proof, which was given in
\cite[Proposition 1.2]{IO}.

A simple argument (a version of Frobenius' lemma) shows that
$$
\Phi(u-\tfrac12;\la)=\prod_{i=1}^{\de(\la)}\frac{u+b_i}{u-a_i}
$$
(this identity can also be deduced from formula \eqref{2C} below). It follows
$$
\log\Phi(u-\tfrac12;\la)=\sum_{m=1}^\infty
\frac{u^{-m}}m\sum_{i=1}^{\de(\la)}\left(a_i^m-(-b_i)^m\right),
$$
which implies that $\A$ is freely generated by the functions
\begin{equation}\label{2B}
p_m(\la):=\sum_{i=1}^{\de(\la)}\left(a_i^m-(-b_i)^m\right), \qquad m=1,2,\dots,
\end{equation}
which are {\it super-power sums\/} in $\{a_i\}$ and $\{-b_i\}$.
\end{proof}

Another characterization of regular functions is provided by

\begin{proposition}\label{prop2.D}
$\A$ coincides with the unital algebra generated by the function $\la\mapsto
|\la|$ and the functions $G_\varphi(\la)$ of the form\/ \eqref{1B}.
\end{proposition}

\begin{proof}
This result is due to S.~Kerov. It is pointed out in his note \cite{Ke1}, see
also \cite[proof of Theorem 4]{KeO}. Here is a detailed proof taken  from
Kerov's unpublished work notes:

We claim that the algebra $\A$ is freely generated by the functions
$$
\widehat p_r(\la)=\sum_{\square\in\la}(c(\square))^r, \qquad r=0,1,\dots\,,
$$
where the sum is taken over the boxes $\square$ of $\la$ and $c(\square)$
denotes the content of a box. Note that $\widehat p_0(\la)=|\la|$.

Indeed, we start with the relation
\begin{equation}\label{2C}
\Phi(u-\tfrac12;\la)=\prod_{i=1}^{\ell(\la)}\frac{u+i-\tfrac12}{u-\la_i+i-\tfrac12}
=\prod_{\square\in\la}\frac{u-c(\square)+\tfrac12}{u-c(\square)-\tfrac12}\,.
\end{equation}
It implies
$$
\log\Phi(u-\tfrac12;\la) =\sum_{m=1}^\infty\frac{u^{-m}}m\sum_{\square\in\la}
\left((c(\square)+\tfrac12)^m-(c(\square)-\tfrac12)^m\right),
$$
or
$$
p_m(\la) =\sum_{k=0}^{\left[\frac{m-1}2\right]}2^{-2k}\binom{m}{2k+1}\widehat
p_{m-1-2k}(\la), \quad m=1,2,\dots,
$$
and our claim follows.
\end{proof}

\begin{remark}
Note a shift of degree: as seen from the above computation, the degree of
$\widehat p_r(\la)$ with respect to the filtration of $\A$ equals $r+1$.
\end{remark}

\begin{remark}\label{rem2.E}
Proposition \ref{prop2.C} makes it possible to introduce a natural algebra
isomorphism between $\La$ and $\A$, which sends the power-sums $p_m\in\La$ to
the functions $p_m(\la)$ defined in \eqref{2B},
\end{remark}

\begin{remark}
The algebra $\A$ is stable under the change of the argument $\la\mapsto\la'$
(transposition of diagrams): this claim is not obvious from the initial
definition but becomes clear from Proposition \ref{prop2.C} or Proposition
\ref{prop2.D}.
\end{remark}

Finally, note that one more characterization of the algebra $\A$ is given in
Section \ref{6}.

\section{A proof of Theorem \ref{thm1.A}}\label{3}

The {\it Young graph\/} has $\Y$ as the vertex set, and the edges are formed by
couples of diagrams that differ by a single box. This is a graded graph: its
$n$th level ($n=0,1,\dots$) is the subset $\Y_n\subset\Y$. The notation
$\mu\nearrow\la$ or, equivalently, $\la\searrow\mu$ means that $\la$ is
obtained from $\mu$ by adding a box (so that the couple $\{\mu,\la\}$ forms an
edge). The quantity $\dim\la$ coincides with the number of monotone paths
$\varnothing\nearrow\dots\nearrow\la$ in the Young graph.

More generally, for any two diagrams $\mu,\la\in\Y$ we denote by
$\dim(\mu,\la)$ the number of monotone paths $\mu\nearrow\dots\nearrow\la$ in
the Young graph that start at $\mu$ and end at $\la$. If there is no such path,
then we set $\dim(\mu,\la)=0$. Equivalently, $\dim(\mu,\la)$ is the number of
standard tableaux of skew shape $\la/\mu$ when $\mu\subseteq\la$, and
$\dim(\mu,\la)=0$ otherwise.

Let $x^{\down m}$ stand for the $m$th falling factorial power of $x$. That is,
$$
x^{\down m}=x(x-1)\dots(x-m+1), \qquad m=0,1,\dots\,.
$$
With an arbitrary $\mu\in\Y$ we associate the following function on $\Y$:
\begin{equation}\label{3A}
F_\mu(\la)=n^{\down m}\, \frac{\dim(\mu,\la)}{\dim \la}\,, \qquad \la\in\Y,
\quad n=|\la|, \quad m=|\mu|.
\end{equation}

\begin{proposition}\label{prop3.A}
For any $\mu\in\Y$, the function $F_\mu$ belongs to $\A$ and has degree
$|\mu|$. Under the isomorphism $\operatorname{gr}\A\simeq\La$, the top degree
term of $F_\mu$ coincides with the Schur function $s_\mu$.
\end{proposition}

\begin{proof}
This can be deduced from \cite[Theorem 5]{KeO}. For direct proofs, see
\cite[Theorem 8.1]{OO1} and \cite[Proposition 1.2]{ORV2}.
\end{proof}

\begin{remark}
Under the isomorphism between $\A$ and $\La^*$, $F_\mu$ turns into the {\it
shifted Schur function\/} $s^*_\mu$, see \cite[Definition 1.4]{OO1}. Under the
isomorphism between $\A$ and $\La$ (Remark \ref{rem2.E}), $F_\mu$ is identified
with the {\it Frobenius--Schur function\/} $Fs_\mu$, see \cite{ORV1},
\cite[Section 2]{ORV2}.
\end{remark}

Introduce a notation for the $n$th Plancherel measure:
\begin{equation}\label{3B}
 M_n(\la)=\frac{(\dim\la)^2}{n!}\,, \qquad \la\in\Y_n\,.
\end{equation}
Thus, the $n$th Plancherel average of a function $F$ on $\Y$ is
\begin{equation}\label{3C1}
\langle F\rangle_n=\sum_{\la\in\Y_n}F(\la)M_n(\la).
\end{equation}

By virtue of Proposition \ref{prop2.D}, Theorem \ref{thm1.A} follows from

\begin{theorem}\label{thm3.B} For any $F\in\A$, $\langle F\rangle_n$ is a polynomial in $n$
of degree at most\/ $\deg F$, where $\deg$ refers to degree with respect to the
filtration in $\A$. Furthermore,
\begin{equation}\label{3C2}
\langle F_\mu\rangle_n=\binom{n}{m}\dim\mu, \qquad \mu\in\Y, \quad m:=|\mu|.
\end{equation}
\end{theorem}

\begin{proof}
First, let us check \eqref{3C2}. If $n<m$ then the both sides of \eqref{3C2}
vanish: the restriction of $F_\mu$ to $\Y_n$ is identically $0$ and
$\binom{n}{m}=0$. Consequently, we may assume $n\ge m$.

Let $(\,\cdot\,,\,\cdot\,)$ denote the standard inner product in $\La$. The
simplest case of Pieri's rule for the Schur functions says that
$$
p_1s_\mu=\sum_{\mu^\bullet:\,\mu^\bullet\searrow\mu}s_{\mu^\bullet}\,.
$$
It follows that for $\la\in\Y_n$
\begin{equation}\label{3D}
\dim(\mu,\la)=(p_1^{n-m}s_\mu,s_\la), \qquad \dim\la=(p_1^n,s_\la).
\end{equation}

Therefore, using the definition \eqref{3A}, we have
\begin{gather*}
\langle F_\mu\rangle_n =\frac{n^{\down m}}{n!}\,\sum_{\la\in\Y_n}
\dim(\mu,\la)\dim\la\\ =\frac{n^{\down m}}{n!}\,\sum_{\la\in\Y_n}
(p_1^{n-m}s_\mu,s_\la)(p_1^n,s_\la) =\frac{n^{\down
m}}{n!}\,(p_1^{n-m}s_\mu,p_1^n)\\ =\frac{n^{\down m}}{n!}\,\left(s_\mu,
\frac{\pd^{n-m}}{\pd p_1^{n-m}}p_1^n\right) =\frac{n^{\down m}}{m!}\,(s_\mu,
p_1^m)=\binom{n}{m}\dim\mu,
\end{gather*}
as required.

By virtue of Proposition \ref{prop3.A}, $\deg F_\mu=|\mu|$ and  $\{F_\mu\}$ is
a basis in $\A$ compatible with the filtration. On the other hand,
$\binom{n}{m}$ is a polynomial in $n$ of degree $m$. Therefore, the first claim
of the theorem follows from \eqref{3C2}.
\end{proof}

\begin{remark}
Stanley \cite[Section 3]{St3} shows that the claim of Theorem \ref{thm1.A}
generalizes to functions of the form $G_\varphi H_\psi$, where $\psi$ is an
arbitrary symmetric function and
\begin{equation}\label{3E}
H_\psi(\la):=\psi(\la_1+|\la|-1,\la_2+|\la|-2,\dots,\la_{|\la|},0,0,\dots),
\qquad \la\in\Y.
\end{equation}
This apparently stronger result also follows from Theorem \ref{thm3.B}, because
(as is readily seen) any function of the form \eqref{3E} belongs to the algebra
$\A$.
\end{remark}

\section{The Jack deformation of the algebra $\A$}\label{4}

Here we extend the definitions of Section \ref{2} by introducing the
deformation parameter $\tth>0$. The previous picture corresponds to the
particular value $\tth=1$. We call $\tth$ the {\it Jack parameter\/}, because
of a close relation to Jack symmetric functions. Note that $\tth$ is inverse to
the parameter $\al$ used in Macdonald's book \cite{Ma} and Stanley's paper
\cite{St1}.

\begin{definition}\label{defn4.A}
The {\it $\tth$-characteristic function\/} of a diagram $\la\in\Y$ is defined
as
$$
\Phi_\tth(u;\la)=\prod_{i=1}^\infty\frac{u+\tth i}{u-\la_i+\tth i}
=\prod_{i=1}^{\ell(\la)}\frac{u+\tth i}{u-\la_i+\tth i}\,.
$$
\end{definition}

This is again a rational function in $u$, regular at infinity and hence
admitting  the Taylor expansion at $u=\infty$ with respect to $u^{-1}$.

\begin{definition}\label{defn4.B}
The algebra $\A_\tth$ of {\it $\tth$-regular functions\/} on $\Y$ is the unital
$\R$-algebra generated by the coefficients of the Taylor expansion at
$u=\infty$ of the function $\Phi_\tth(u;\la)$ (or, equivalently, of
$\log\Phi_\tth(u;\la)$).
\end{definition}

The Taylor expansion of $\log\Phi_\tth(u;\la)$ at $u=\infty$ has the form
$$
\log\Phi_\tth(u;\la)=\sum_{m=1}^\infty\frac{p^*_{m;\tth}(\la)}m\,u^{-m},
$$
where, by definition,
$$
p^*_{m;\tth}(\la)=\sum_{i=1}^\infty [(\la_i-\tth i)^m-(-\tth i)^m], \qquad
m=1,2,\dots, \quad \la\in\Y
$$
(as above, summation actually can be taken up to $i=\ell(\la)$). Thus, the
algebra $\A_\tth$ is generated by the functions
$p^*_{1;\tth},p^*_{2;\tth},\dots$. These functions are algebraically
independent.

The filtration in $\A_\tth$ is introduced exactly as in the particular case
$\tth=1$. We still have a canonical isomorphism of graded algebras
$\gr(\A_\tth)\simeq\La$ and a canonical isomorphism of filtered algebras
$\A\simeq\La^*_\tth$, where $\La^*_\tth$ denotes the algebra of $\tth$-shifted
symmetric functions \cite{KOO}. However, for general $\tth$, we do not see a
natural way to define an isomorphism between $\A_\tth$ and $\La$.

\section{Jack deformation of Plancherel averages}\label{5}

Recall that $\tth>0$ is a fixed parameter, which is inverse to Macdonald's
\cite{Ma} parameter $\al$. We consider the Jack deformation
$(\,\cdot\,,\,\cdot\,)_\tth$ of the standard inner product in the algebra $\La$
of symmetric functions. In the basis $\{p_\la\}$ of power-sum functions,
\begin{equation}\label{5A}
(p_\la,p_\mu)_\tth=\delta_{\la\mu}z_\la\tth^{-|\la|}, \qquad \la,\mu\in\Y,
\end{equation}
cf. \cite[Chapter VI, Section 10]{Ma}; the standard notation $z_\la$ is
explained in \cite[Chapter I, Section 2]{Ma}. Let $\{P_\la\}$ and $\{Q_\la\}$
be the biorthogonal bases formed  the $P$ and $Q$ Jack symmetric functions
(which differ from each other by normalization factors). In Macdonald's
notation (\cite[Chapter VI, Section 10]{Ma}), these are $P^{(1/\tth)}_\la$ and
$Q^{(1/\tth)}_\la$. To simplify the notation, we will not include $\tth$ into
the notation for the Jack functions. When $\tth=1$, the both versions of the
Jack functions turn into the Schur functions $s_\la$.

Introduce the notation
\begin{equation}\label{5B}
\dim_\tth \la=(p_1^n,Q_\la)_\tth, \quad \dim'_\tth \la=(p_1^n,P_\la)_\tth\,,
\qquad \la\in\Y_n\,.
\end{equation}
More generally, we set (cf. \eqref{3D})
\begin{equation}\label{5C1}
\dim_\tth(\mu,\la)=(p_1^{|\la|-|\mu|}P_\mu,Q_\la)_\tth, \quad
\dim'_\tth(\mu,\la)=(p_1^{|\la|-|\mu|}Q_\mu,P_\la)_\tth,
\end{equation}
where we assume $|\mu|\le|\la|$; otherwise the dimension is set to be 0.

\begin{proposition}\label{prop5.A} The quantities \eqref{5B} are strictly positive.
The quantities \eqref{5C1} are strictly positive if $\mu\subseteq\la$ and
vanish otherwise.
\end{proposition}

\begin{proof}
The first claim being a particular case of the second one, we focus on the
second claim. We employ the formalism described in \cite{KOO}.

The simplest case of Pieri's rule for Jack symmetric functions (\cite[Chapter
VI, Section 10 and (6.24)(iv)]{Ma}) says that $p_1P_\mu$ is a linear
combination of the functions $P_{\mu^\bullet}$, $\mu^\bullet\searrow\mu$, with
strictly positive coefficients. The coefficients are just the quantities
$\ka_\tth(\mu,\mu^\bullet):=(p_1P_\mu,Q_{\mu^\bullet})_\tth$; let us view them
as {\it formal multiplicities\/} attached to the edges
$\mu\nearrow\mu^\bullet$. More generally, the {\it weight\/} of a finite
monotone path $\mu\nearrow\dots\nearrow\la$ in the Young graph is defined as
the product of the formal multiplicities of edges entering the path. Observe
now that $\dim_\tth(\mu,\la)$ is the sum of the weights of all monotone paths
connecting $\mu$ to $\la$. This proves the claim concerning
$\dim_\tth(\mu,\la)$. For $\dim'_\tth(\mu,\la)$ the argument is the same: we
simply swap the $P$ and $Q$ functions.
\end{proof}

With an arbitrary $\mu\in\Y$ we associate the following function on $\Y$, cf.
\eqref{3A}:
\begin{equation*}
F_{\mu;\tth}(\la)=n^{\down m}\, \frac{\dim_\tth(\mu,\la)}{\dim_\tth \la}\,,
\qquad \la\in\Y, \quad n=|\la|, \quad m=|\mu|.
\end{equation*}

\begin{proposition}\label{prop5.B}
For any $\mu\in\Y$, the function $F_{\mu;\tth}$ just defined belongs to
$\A_\tth$. Under the isomorphism $\operatorname{gr}\A_\tth\simeq\La$, the top
degree term of $F_{\mu;\tth}$ coincides with the Jack function $P_\mu$.
\end{proposition}

\begin{proof} See \cite[Section 5]{OO2}. Note that under the isomorphism
$\La_\tth^*\to\A_\tth$, $F_{\mu;\tth}$ coincides with the image of the {\it
shifted Jack function\/} $P^*_\mu$.
\end{proof}

\begin{definition}
The {\it Jack deformation of the Plancherel measure\/} with parameter $\tth$ on
the set $\Y_n$ (or {\it Jack--Plancherel measure\/}, for short) is defined by
\begin{equation}\label{5F}
M_{n;\tth}(\la)=\frac{(p_1^n,Q_\la)_\tth(p_1^n,P_\la)_\tth}{(p_1^n,p_1^n)_\tth}\,,
\qquad \la\in\Y_n\,.
\end{equation}
By Proposition \ref{prop5.A}, the quantity $M_{n;\tth}(\la)$ is always
positive. Since $\{P_\la\}$ and $\{Q_\la\}$ are biorthogonal bases, the sum of
the quantities \eqref{5F} over $\la\in\Y_n$ equals 1. Therefore, $M_{n;\tth}$
is a probability measure. Note that the above definition agrees with that given
in \cite[Section 7]{Ke2} and \cite[Section 3.3.2]{Ok}.
\end{definition}

Because
$$
(p_1^n,p_1^n)_\tth=z_{(1^n)}\tth^{-n}=\frac{n!}{\tth^n}\,,
$$
\eqref{5F} can be rewritten as
\begin{equation}\label{5E}
M_{n;\tth}(\la)=\frac{\tth^n(p_1^n,Q_\la)_\tth(p_1^n,P_\la)_\tth}{n!}
=\frac{\tth^n\,\dim_\tth\la\,\dim'_\tth\la}{n!}\,, \quad \la\in\Y_n\,.
\end{equation}

Clearly, for $\tth=1$ the definition coincides with \eqref{3B}.

\begin{remark}
{}From the Jack version of the duality map $\La\to\La$ (\cite[Chapter VI,
(10.17)]{Ma}) it can be seen that under the involution $\la\mapsto\la'$ the
measure $M_{n;\tth}$ is transformed into $M_{n;\tth^{-1}}$.
\end{remark}

Given a function $F$ on $\Y$, its $n$th {\it Jack--Plancherel average\/} is
defined by analogy with \eqref{3C1}:
\begin{equation}\label{5C2}
\langle F\rangle_{n;\tth}=\sum_{\la\in\Y_n}F(\la)M_{n;\tth}(\la).
\end{equation}

Here is a generalization of Theorem \ref{thm3.B}:

\begin{theorem}\label{thm5.C}
For any $F\in\A_\tth$, $\langle F\rangle_{n;\tth}$ is a polynomial in $n$ of
degree at most\/ $\deg F$, where $\deg$ refers to degree with respect to the
filtration in $\A_\tth$. Furthermore,
\begin{equation*}
\langle F_{\mu;\tth}\rangle_{n;\tth}=\tth^m\binom{n}{m}\dim_\tth\mu.
\end{equation*}
\end{theorem}

\begin{proof}
The argument relies on Proposition \ref{prop5.B} and is the same as in the
proof of Theorem \ref{3B}, with minor obvious modifications. In particular, we
use the fact that the adjoint to multiplication by $p_1$ is equal to
$\tth^{-1}\pd/\pd p_1$. For reader's convenience, we repeat the main
computation:
\begin{gather*}
\langle F_{\mu;\tth}\rangle_{n;\tth} =\tth^n\frac{n^{\down
m}}{n!}\,\sum_{\la\in\Y_n} \dim_\tth(\mu,\la)\dim'_\tth\la\\
=\tth^n\frac{n^{\down m}}{n!}\,\sum_{\la\in\Y_n}
(p_1^{n-m}P_\mu,Q_\la)_\tth(p_1^n,P_\la)_\tth =\tth^n\frac{n^{\down
m}}{n!}\,(p_1^{n-m}P_\mu,p_1^n)_\tth\\ =\tth^n\frac{n^{\down m}}{n!}\,(P_\mu,
(\tth^{-1}\pd/\pd p_1)^{n-m}p_1^n)_\tth =\tth^m\frac{n^{\down m}}{m!}\,(P_\mu,
p_1^m)_\tth=\tth^m\binom{n}{m}\dim_\tth\mu\,.
\end{gather*}
\end{proof}

\section{Kerov's interlacing coordinates}\label{6}

Let $\la\in\Y$ be a Young diagram drawn according to the ``English picture''
\cite[Chapter I, Section 1]{Ma}, that is, the first coordinate axis (the row
axis) is directed downwards and the second coordinate axis (the column axis) is
directed to the right. Consider the border line of $\la$ as the directed path
coming from $+\infty$ along the second (horizontal) axis, next turning several
times alternately down and to the left, and finally going away to $+\infty$
along the first (vertical) axis. The corner points on this path are of two
types: the {\it inner corners\/}, where the path switches from the horizontal
direction to the vertical one, and the {\it outer corners\/} where the
direction is switched from vertical to horizontal. Observe that the inner and
outer corners always interlace and the number of inner corners always exceeds
by 1 that of outer corners. Let $2d-1$ be the total number of the corners and
$(r_i,s_i)$, $1\le i\le 2d-1$, be their coordinates. Here the odd and even
indices $i$ refer to the inner and outer corners, respectively.

\medskip

\begin{center}
 \scalebox{0.5}{\includegraphics{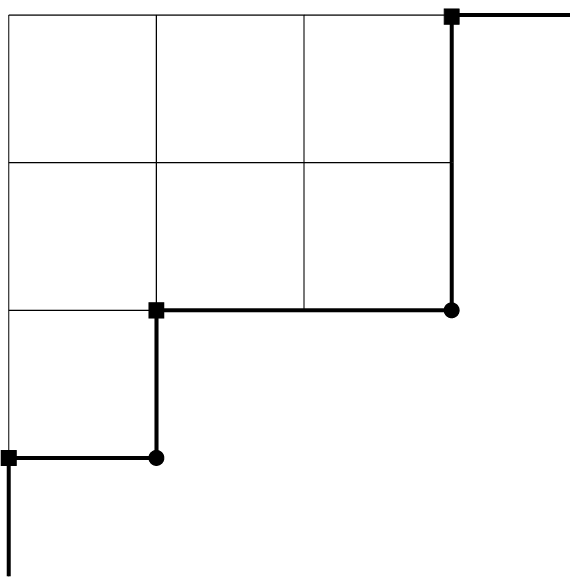}}

Figure 1. The corners of the diagram $\la=(3,3,1)$.
\end{center}

For instance, the diagram $\la=(3,3,1)$ shown on the figure has $d=3$, three
inner corners $(r_1,s_1)=(0,3)$, $(r_3,s_3)=(2,1)$, $(r_5,s_5)=(3,0)$, and two
outer corners $(r_2,s_2)=(2,3)$, $(r_4,s_4)=(3,1)$.

As above, $\tth$ is assumed to be a fixed strictly positive parameter. The
numbers
\begin{equation}\label{6A}
 x_1:=s_1-\tth  r_1, \quad y_1:=s_2-\tth r_2,\,\dots,\,
 y_{d-1}:=s_{2d-2}-\tth r_{2d-2}, \quad x_d:=s_{2d-1}-\tth r_{2d-1}
\end{equation}
form two interlacing sequences of integers
$$
x_1>y_1>x_2>\dots>y_{d-1}>x_d
$$
satisfying the relation
\begin{equation}\label{6B}
\sum_{i=1}^d x_i-\sum_{j=1}^{d-1}y_j=0.
\end{equation}

For instance, if $\la=(3,3,1)$ as in the example above, then
$$
x_1=3, \quad y_1=3-2\tth, \quad x_2=1-2\tth, \quad y_2=1-3\tth, \quad
x_3=-3\tth.
$$

\begin{definition}\label{defn6.A} The two interlacing sequences
\begin{equation}
X=(x_1,\dots,x_d), \quad Y=(y_1,\dots,y_{d-1})
\end{equation}
as defined above are called the ($\tth$-dependent) {\it Kerov interlacing
coordinates\/} of a Young diagram $\la$. (Note that in the case $\tth=1$,
Kerov's $(X,Y)$ coordinates are similar to Stanley's ``$(p,q)$ coordinates''
introduced in \cite{St2}: the two coordinate systems are related by a simple
linear transformation.)
\end{definition}

Let $u$ be a complex variable. Given a Young diagram $\la$, we set
$$
\HH(u;\la)=\frac{u\prod\limits_{j=1}^{d-1}(u-y_j)}{\prod\limits_{i=1}^d(u-x_i)}\,,
$$
and
$$
\p_m(\la)=\sum_{i=1}^d x_i^m-\sum_{j=1}^{d-1}y_j^m, \qquad m=1,2,\dots,
$$
where $X=\{x_i\}$, $Y=\{y_j\}$, and $d$ are as in Definition \ref{defn6.A}.
Obviously,
$$
\log\HH(u;\la)=\sum_{m=1}^\infty \frac{\p_m(\la)}m\,u^{-m}.
$$
Note that $\p_1(\la)\equiv0$ because of \eqref{6B}.

\begin{proposition}\label{prop6.D} The following relation holds
$$
\HH(u;\la)=\frac{\Phi(u-\tth;\la)}{\Phi(u;\la)}\,.
$$
\end{proposition}

\begin{proof} See \cite[Proposition 6.3]{Ol}.
\end{proof}

From this result one deduces:

\begin{proposition}
The functions\/ $\p_m(\la)$ belong to the algebra $\A_\tth$. More precisely, we
have
$$
\p_m=\tth\cdot m\cdot p^*_{m-1;\tth}\,+\,\dots, \qquad m=2,3,\dots,
$$
where dots stand for lower degree terms, which are a linear combination of
elements $p^*_{l;\tth}$ with $1\le l\le m-2$.
\end{proposition}

\begin{proof} See \cite[Proposition 6.5]{Ol}.
\end{proof}

\begin{corollary}  The  functions
$\{\p_2,\p_3,\dots\}$ form a system of algebraically independent generators of
the algebra $\A_\tth$, compatible with the filtration. More precisely, under
the identification of\/ $\A_\tth$ with the algebra of polynomials
$\R[\p_2,\p_3,\dots]$, the filtration is determined by setting
$$
\deg \p_m=m-1, \qquad m=2,3,\dots\,.
$$
\end{corollary}

Thus, the algebra $\A_\tth$ of $\tth$-regular functions coincides with the
algebra of {\it super-symmetric functions\/} in Kerov's $\tth$-dependent
interlacing coordinates.

Consider the expansion in partial fractions for $u^{-1}\HH(u;\la)$:
$$
\frac{\prod\limits_{j=1}^{d-1}(u-y_j)}{\prod\limits_{i=1}^d(u-x_i)}=\sum_{i=1}
^d\frac { \pi^\up_i } { u-x_i}\,.
$$
Here the coefficients $\pi^\up_i$ are given by the formula
$$
\pi^\up_i=\pi^\up_i(\la)=\frac{\prod\limits_{j=1}^{d-1}(x_i-y_j)}
{\prod\limits_{l:\, l\ne i}(x_i-x_l)}\,, \qquad i=1,\dots,d.
$$

Observe that the boxes that may be appended to $\la$ are associated, in a
natural way, with the inner corners of the boundary of $\la$. Consequently, we
may also associate these boxes with the $x$'s: $\square_i\leftrightarrow x_i$.

It is ready to check that the coefficients $\pi^\up_i$ are strictly positive
and sum up to 1. Introduce the notation
$$
p^\up_{n;\tth}(\la,\la\cup\square_i)=\pi_i^\up(\la), \qquad 1\le i\le d, \quad
\la \in\Y_n
$$
(the quantities $\pi^\up_i(\la)$ in the right-hand side depend on $\tth$
through \eqref{6A}). We regard $p^\up_{n;\tth}$ as a transition function acting
from $\Y_n$ to $\Y_{n+1}$. The system $\{p^\up_{n;\tth}\}_{n=0,1,\dots}$
determines a model of random growth of Young diagrams: an inhomogeneous Markov
chain on $\Y$ whose state at time $n=0,1,\dots$ is a diagram from $\Y_n$. Every
trajectory of this Markov chain is an infinite monotone path in $\Y$ starting
at $\varnothing$.

Denote by $M'_{n;\tth}$ the marginal distribution of this Markov chain after
$n$ steps. That is, $M'_{n;\tth}$ is the probability measure on $\Y_n$ defined
by the recursion
\begin{equation}\label{6C}
M'_{n+1;\tth}(\nu)=\sum_{\la\in\Y_n:\,\la\nearrow\nu}M'_{n;\tth}(\la)p^\up_{n;\tth}(\la,\nu)
\end{equation}
with the initial condition $M'_{0;\tth}(\varnothing)=1$.

\begin{proposition}\label{prop6.B}  $M'_{n;\tth}$ coincides with the Jack--Plancherel
measure $M_{n;\tth}$ as defined in \eqref{5F}
\end{proposition}

\begin{proof}
This is one of the main results of Kerov \cite{Ke2} (see Section 7 in
\cite{Ke2}). For $\tth=1$, it allows a direct elementary verification. For
general $\tth$, the proof given in \cite{Ke2} is more delicate; it uses the
hook-type formulas for $\dim_\tth\la$ and $\dim'_\tth\la$ (see \cite[Section
6]{Ke2} and \cite[Section 5]{St1}).
\end{proof}

If we agree to take \eqref{6C} as the initial definition of the Jack
deformation of the Plancherel measure, then (as will be seen) we may completely
eliminate the Jack polynomials from our considerations.

Let us restate the first claim of Theorem \ref{thm5.C} in terms of the measure
$M'_{n;\tth}$:

\begin{theorem}\label{thm6.C}
Let\/ $\langle\,\cdot\,\rangle'_{n;\tth}$ stand for the expectation with
respect to the measure $M'_{n;\tth}$. For any $F\in\A_\tth$, $\langle
F\rangle'_{n;\tth}$ is a polynomial in $n$ of degree at most\/ $\deg F$.
\end{theorem}

We will deduce Theorem \ref{thm6.C} from the following claim.

Let $\pd$ denote the operator acting in the space of functions on $\Y$ as
\begin{equation}\label{6D}
(\pd F)(\la)=-F(\la)+\sum_{\nu:\, \nu\searrow\la}
p^\up_{n;\tth}(\la,\nu)F(\nu), \qquad \la\in\Y, \quad n=|\la|.
\end{equation}

\begin{theorem}\label{thm6.D} The operator $\pd$ defined by \eqref{6D} preserves
the algebra $\A_\tth$ and reduces degree by\/ $1$.
\end{theorem}

\begin{proof}[Reduction of Theorem \ref{thm6.C} to Theorem \ref{thm6.D}] By
virtue of \eqref{6C},
$$
\langle F\rangle'_{n+1;\tth}-\langle F\rangle'_{n;\tth}=\langle \pd
F\rangle'_{n;\tth}\,, \qquad n=0,1,\dots\,.
$$
Since $\langle 1\rangle'_{n;\tth}=1$, the claim of Theorem \ref{thm6.C} is
obtained by induction on $\deg F$.
\end{proof}

\begin{proof}[Proof of Theorem \ref{thm6.D}]
The claim of Theorem \ref{thm6.D} can be obtained by a degeneration from a much
more general claim, \cite[Theorem 7.1(ii)]{Ol}. In the notation of \cite{Ol},
the degeneration consists in letting  certain parameters $z$ and $z'$ go to
infinity.

An alternative possibility is to adapt the approach of \cite{Ol} to the present
situation by eliminating these parameters at all, which substantially
simplifies the computations. Here is a sketch of the argument; for more detail
we refer to \cite{Ol}.

Introduce the functions $\h_0(\la), \h_1(\la), \dots$ on $\Y$ from the
decomposition
\begin{equation*}
\HH(u;\la)=\sum_{m=0}^\infty \h_m(\la)u^{-m}
\end{equation*}
and note that
\begin{equation*}
\h_0(\la)\equiv1, \quad \h_1(\la)\equiv0.
\end{equation*}
The functions $\h_2,\h_3,\dots$ are algebraically independent generators of the
algebra $\A_\tth$.

For a partition $\rho=(\rho_1,\rho_2,\dots)$, set
\begin{equation*}
\h_\rho=\h_{\rho_1}\h_{\rho_2}\dots\,.
\end{equation*}
Because of $\h_1=0$, we will assume in what follows that $\rho$ does not have
parts $\rho_i$ equal to 1 (otherwise $\h_\rho=0$). Then the elements $\h_\rho$
form a linear basis in $\A_\tth$, consistent with filtration:
\begin{equation*}
\deg \h_\rho=|\rho|-\ell(\rho),
\end{equation*}
where $|\rho|=\sum\rho_i$ and $\ell(\rho)$ is the number of nonzero parts in
$\rho$. This is related to the fact that $\deg\h_m=m-1$ for $m=2,3,\dots$\,.

We aim at computing the action of $\pd$ on the basis elements $\h_\rho$. For
any $k=1,2,\dots$, we have
\begin{equation*}
\prod_{l=1}^k \HH(u_l;\la) =\sum_{\rho:\, \ell(\rho)\le
k}m_\rho(u_1^{-1},\dots,u_k^{-1})\h_\rho(\la),
\end{equation*}
where $m_\rho$ is the monomial symmetric function. Thus, we may view finite
products $\prod_{l}\HH(u_l;\la)$ as generating series for the basis elements
$\h_\rho$. It is convenient to consider first the action of $\pd$ on these
generating series.

The argument in \cite[Section 7.2]{Ol} shows that $1+\pd$ acts on
$\HH(u_1;\la)\dots\HH(u_k;\la)$ as multiplication by the series
\begin{equation*}
F^\up(u_1,\dots,u_k;\la):=\sum_{i=1}^d \pi^\up_i(\la)
\prod_{l=1}^k\frac{(u_l-x_i)(u_l-x_i+\tth-1)}{(u_l-x_i-1)(u_l-x_i+\tth)}\,.
\end{equation*}
This series belongs to $\A_\tth[[u_1^{-1},\dots,u_l^{-1}]]$, because of the
fundamental identity
\begin{equation*}
\sum_{i=1}^d \pi^\up_i(\la)x_i^m=\h_m(\la), \quad m=0,1,\dots,
\end{equation*}
see \cite[Lemma 6.11]{Ol}. It follows that $\pd$ maps $\A_\tth$ into itself.

A more detailed analysis (see \cite[Section 7.3]{Ol}) shows the following.
Introduce the linear map $f\to\langle f\rangle^\up$ from $\R[x]$ to $\A_\tth$
by setting
\begin{equation*}
\langle x^m\rangle^\up=\h_m.
\end{equation*}
Next, write the decomposition
\begin{equation*}
\prod_{l=1}^k\frac{(u_l-x)(u_l-x+\tth-1)}{(u_l-x-1)(u_l-x+\tth)} =\sum_{\si}
a_\si(x)m_\sigma(u_1^{-1},\dots,u_k^{-1}),
\end{equation*}
where $\si$ ranges over partitions with $\ell(\si)\le k$ and $a_\si(x)$ are
appropriate polynomials. Finally, let $c^\rho_{\si\tau}$ be the structure
constants of the algebra $\La$ in the basis of monomial symmetric functions:
$$
m_\si m_\tau=\sum_\rho c^\rho_{\si\tau}m_\rho
$$
(note that $c^\rho_{\si\tau}$ vanishes unless $|\rho|=|\si|+|\tau|$). Then we
have
\begin{equation}\label{6E}
 (1+\pd)\h_\rho=\sum_{\si,\tau:\,
|\si|+|\tau|=|\rho|}c^\rho_{\si\tau}\left\langle
a_\si(x)\right\rangle^\up\h_\tau,
\end{equation}
see \cite[Lemma 7.4]{Ol}. Note that
\begin{equation}\label{6F}
a_\si(x)=a_{\si_1}(x)a_{\si_2}(x)\dots,
\end{equation}
where
\begin{gather}
a_s(x)=(s-1)\tth x^{s-2}+\frac{(s-1)(s-2)}2\,\tth(1-\tth)x^{s-3}\,+\dots,
\qquad s\ge 2  \label{6G}\\
\langle a_0(x)\rangle^\up\equiv1, \quad \langle a_1(x)\rangle^\up\equiv0
\label{6H}
\end{gather}
(see \cite[Lemma 7.3]{Ol}). It follows, in particular, that we may assume that
in \eqref{6E}, $\si$ does not have parts equal to 1.

Identify $\A_\tth$ with the polynomial algebra $\R[\h_2,\h_3,\dots]$. Using the
argument of \cite[Lemma 7.12]{Ol}, one deduces from formulas \eqref{6E},
\eqref{6F}, \eqref{6G}, and \eqref{6H} that
\begin{equation*}
\pd=\tth\frac{\pd}{\pd\h_2}+\dots
\end{equation*}
where the dots stand for terms of degree $\le-2$ (that is, operators in
$\A_\tth$ reducing degree at least by 2). This concludes the proof, since the
operator $\pd/\pd\h_2$ reduces degree by 1 (recall that $\h_2$ has degree 1).
\end{proof}


\begin{thebibliography}{XXXX}

\bibitem{Fuj}
S. Fujii, H. Kanno, S. Moriyama and S. Okada, {\it Instanton calculus and
chiral one-point functions in supersymmetric gauge theories\/}. Advances Theor.
Math. Phys. {\bf12} (2008) 1401--1428.

\bibitem{Han}
G.-N. Han, {\it Some conjectures and open problems on partition hook
lengths\/}. Exper. Math. {\bf18} (2009), 97--106.

\bibitem{IO}
V. Ivanov and G. Olshanski, {\it Kerov's central limit theorem for the
Plancherel measure on Young diagrams\/}. In: Symmetric functions 2001. Surveys
of developments and perspectives. Proc. NATO Advanced Study Institute
(S.~Fomin, editor), Kluwer, 2002, pp. 93--151; arXiv:math/0304010.

\bibitem{Ke1}
S. Kerov, {\it Gaussian limit for the Plancherel measure of the symmetric
group\/}. Comptes Rendus Acad. Sci. Paris, S\'er. I {\bf316} (1993), 303--308.

\bibitem{Ke2}
S. Kerov, {\it Anisotropic Young diagrams and Jack symmetric functions\/}.
Function. Anal. i Prilozhen. {\bf34} (2000), no.~1, 51--64 (Russian); English
translation: Funct. Anal. Appl. {\bf34} (2000), 45--51; arXiv:math/9712267.

\bibitem{KOO}
S.~Kerov, A.~Okounkov, and G.~Olshanski,  {\it The boundary of Young graph with
Jack edge multiplicities\/}. Intern. Math. Res. Notices (1998),  no.~4,
173--199; arXiv: q-alg/9703037.

\bibitem{KeO}
S. Kerov and G. Olshanski, {\it Polynomial functions on the set of Young
diagrams\/}. Comptes Rendus Acad. Sci. Paris, Ser. I {\bf319} (1994), 121--126.

\bibitem{Ma}
I.~G.~Macdonald, {\it Symmetric functions and Hall
polynomials\/}, 2nd edition.  Oxford University Press, 1995.

\bibitem{Ok}
A. Okounkov, {\it The uses of random partitions\/}. In:
XIVth International Congress on Mathematical Physics,   World Sci. Publ.,
Hackensack, NJ, 2005, pp. 379--403; arXiv:math-ph/0309015.

\bibitem{OO1}
A. Okounkov and G. Olshanski, {\it Shifted Schur functions\/}, Algebra i Analiz
{\bf9} (1997), no. 2, 73--146 (Russian); English version: St. Petersburg
Mathematical J., {\bf9} (1998), 239--300; arXiv:q-alg/9605042.


\bibitem{OO2}
A. Okounkov and G. Olshanski, {\it Shifted Jack polynomials, binomial formula,
and applications\/}. Math. Research Lett. {\bf4} (1997), 69--78;
arXiv:q-alg/9608020.

\bibitem{Ol}
G. Olshanski, {\it Anisotropic Young diagrams and infinite-dimensional
diffusion processes with the Jack parameter\/}. Intern. Math. Research Notices,
to appear; arXiv:0902.3395.

\bibitem{ORV1}
G. Olshanski, A. Regev, and A. Vershik, {\it Frobenius--Schur functions:
summary of results\/}, arXiv:math/0003031.

\bibitem{ORV2}
G. Olshanski, A. Regev, and A. Vershik, {\it Frobenius--Schur functions\/}. In:
Studies in Memory of Issai Schur (A.~Joseph, A.~Melnikov, R.~Rentschler, eds),
Progress in Mathematics {\bf 210}, Birkh\"auser, 2003, pp. 251--300;
arXiv:math/0110077.

\bibitem{St1}
R. P. Stanley, {\it Some combinatorial properties of Jack symmetric
functions\/}. Adv. Math. {\bf 77} (1989), 76--115.

\bibitem{St2}
R. P. Stanley,  {\it Irreducible symmetric group characters of rectangular
shape\/}. S\'emin. Lothar. Combin.  {\bf50} (2003),  Art. B50d, 11 pp.

\bibitem{St3}
R. P. Stanley, {\it Some combinatorial properties of hook lengths, contents,
and parts of partitions\/}. Ramanujan J., to appear;  arXiv:0807.0383.

\bibitem{VK}
A.~M.~Vershik and S.~V.~Kerov, {\it Asymptotic theory of characters of the
symmetric group\/}. Function. Anal. i Prilozhen. {\bf15} (1981), no.~4, 15--27
(Russian); English translation: Funct. Anal. Appl. {\bf 15} (1981), 246--255.

\end{thebibliography}
\end{document}